\documentclass[12pt]{amsart}

\usepackage{amsfonts,amsmath,latexsym,amssymb,verbatim,amsbsy,times, esint}
\usepackage{amsthm}
\usepackage{pstricks}
\usepackage{mathtools}
\usepackage{graphicx}
\usepackage{caption}
\usepackage[normalem]{ulem}

\usepackage[top=0.8 in, bottom=0.8in, left=0.8in, right=0.8in]{geometry}

\usepackage[shortlabels]{enumitem}
\usepackage[colorlinks=true, pdfstartview=FitV, linkcolor=blue,citecolor=green, urlcolor=blue]{hyperref}

\usepackage{enumitem}

\theoremstyle{plain}
\newtheorem{THEOREM}{Theorem}[section]
\newtheorem{COROL}[THEOREM]{Corollary}

\newtheorem{PROP}[THEOREM]{Proposition}

\theoremstyle{definition}
\newtheorem{DEF}[THEOREM]{Definition}

\theoremstyle{remark}
\newtheorem{REMARK}[THEOREM]{Remark}

\newcommand{\R}{\ensuremath{\mathbb{R}}}   

\renewcommand{\S}{\ensuremath{\mathbb{S}}}

\def \d {\delta}
\def \g {\gamma}

\def \f {\varphi}

\def \k {\kappa}
\def \l {\lambda}
\def \L {\Lambda}

\def \n {\nabla}
\def \s {\sigma}

\def \t {\tau}
\def \w {\omega}

\def \bn {{\bf n}}

\def \bu {{\bf u}}

\def\cB {\mathcal {B}}

\def\cF {\mathcal {F}}

\def \cP {\mathcal{P}}

\def \p {\partial}

\DeclareMathOperator{\Span}{span} %
\DeclareMathOperator{\supp}{supp} %
\DeclareMathOperator{\diver}{div} %

\def \cP {\mathcal{P}}

\def \ds  {\, \mathrm{d}s}

\def \dtau  {\, \mathrm{d}\tau}

\def \dx  {\, \mathrm{d}x}

\def \ds  {\, \mathrm{d}s}

\def \sfS  { \mathsf{S}}
\def \sfT  { \mathsf{T}}

\def \brho {\bar{\rho}}

\def \eff {\mathrm{eff}}

	\title[Self-gravitating hyperelastic matter]{Separable motions for self-gravitating hyperelastic matter}
	
	\author{Juhi Jang}
	\address[J. Jang]{Department of Mathematics, University of Southern California, Los Angeles, CA 90089}
	\email{juhijang@usc.edu}
	
	\author{Trevor M. Leslie}
	\address[T. Leslie]{Department of Applied Mathematics, Illinois Institute of Technology, Chicago, IL 60616}
	\email{tleslie@iit.edu}

\begin{document}


\begin{abstract}
	In this paper, we prove the existence of separable solutions to the equations of motion for self-gravitating hyperelastic matter, under an appropriate class of constitutive assumptions on the strain-energy function.  Our framework includes both global-in-time solutions which expand 
	and also solutions which collapse to a point in finite time. Other authors have constructed expanding solutions in similar settings, but to the best of our knowledge, the collapsing solutions we construct are completely 
	new.  
\end{abstract}	

\dedicatory{Dedicated to Professor Thomas C. Sideris on the Occasion of His 70th Birthday}

\maketitle

\section{Introduction} 

In this paper, we consider the dynamics of a self-gravitating elastic  body $\cB(t)$, surrounded by vacuum, which obeys the following equations of motion:
\begin{equation}
	\label{e:Eulerian}
	\begin{cases} 
		\p_t \rho(t,x) + \diver(\rho \bu)(t,x) = 0,
		\qquad x\in \cB(t) = \supp \rho(t), \\
		\rho (\p_t\bu + \bu\cdot \n \bu)(t,x) = \diver \sfT(t,x) - G\rho(t,x) \n_x \Phi[\rho(t)](x), \qquad x\in \cB(t),\\
		\Delta \Phi[\varrho] = 4\pi \varrho, \qquad \lim_{|x|\to \infty} \Phi[\varrho](x) = 0 \\
		\sfT(t,x)\bn(x) = 0, \qquad x\in \p \cB(t).
	\end{cases}
\end{equation}
Here $\rho\ge 0$ and $\bu$ are the density and velocity of the body, $\sfT$ is the Cauchy stress tensor, and $\bn(x)$ denotes the unit normal vector to the boundary of $\cB(t) = \supp \rho(t)$.   The constant $G$ is Newton's gravitational constant.  In order to determine the evolution of the quantities $\rho$ and $\bu$, the system must be supplemented with a constitutive relation that describes how the Cauchy stress $\sfT$ responds to the dynamics; specifying a suitable class of constitutive hypotheses is an important and nontrivial part of formulating and studying problems related to \eqref{e:Eulerian}.  We discuss our choice of constitutive assumptions in Sections \ref{ss:constitutive} and \ref{ss:results}.

Our main goal is to construct a family of special spherically symmetric, `separable' solutions to the system \eqref{e:Eulerian}, under an appropriate set of constitutive assumptions.  Our framework is in large part inspired by a recent paper of Sideris \cite{Sideris2023}, who constructed separable solutions in the absence of a gravitational term.  Calogero \cite{Calogero2022} has also studied separable solutions of \eqref{e:Eulerian}, though there are some important differences between his framework and ours.  We will discuss the relationship between our work and the papers \cite{Sideris2023}, \cite{Calogero2022} in Section \ref{ss:results}, after providing more details about the setup of our problem.  

\subsection{Reference Configuration and Separable Solutions} Our constitutive assumptions and the separable solutions we will construct are most naturally described in Lagrangian rather than Eulerian coordinates.  We fix the following reference domain $\cB$ and (constant) reference density $\brho$:
\begin{equation}
	\label{e:refdomain} 
	\cB = \{X\in \mathbb{R}^3 : |X|\le 1\},
\end{equation} 
\begin{equation} 
	\label{e:refdensity}
	\brho = \frac{1}{|\cB(0)|}\int_{\cB(0)} \rho(0,x) \dx. 
\end{equation}
Here $|\cB(0)|$ denotes Lebesgue measure of $\cB(0)$.  When we refer to a `flow map' associated to a classical solution $(\rho, \bu)$ of \eqref{e:Eulerian} on the time interval $[0,T)$ (assumed to be maximal), we mean a function $x~:~[0,T)\times \cB \to \cB(t)$ that satisfies $\dot{x}(t,X) = \bu(t,x(t,X))$ and $\rho(t, x(t,X)) \det Dx(t,X) = \brho$, for all $t\in [0,T)$ and all $X\in \cB$.  We also require that $x(t, \cdot)$ is an orientation-preserving $C^1$ diffeomorphism for each $t\in [0,T)$ and that $x$ is at least $C^2$ in time.  

Any classical spherically symmetric solution will possess a flow map satisfying $x(t,X) = \phi(t,|X|)\frac{X}{|X|}$ for $t\in [0,T)$ and $X\in \cB\backslash \{0\}$, for some function $\phi:[0,T)\times [0,1]\to [0,+\infty)$.  In order for $x(t, \cdot)$ to be an orientation-preserving, $C^1$ diffeomorphism for each $t$, we require that $\phi(t,0) = \phi_{RR}(t,0) = 0$ for all $t\in [0,T)$ and $\phi_R(t, R)>0$ on $[0,T)\times [0,1]$ (c.f. \cite{Sideris2023}, Section 5).  (Here $\phi_R$ and $\phi_{RR}$ denote partial derivatives of $\phi$.)

When we refer to a \textit{separable solution}, we mean a spherically symmetric solution such that the function $\phi(t,R)$ from the previous paragraph has the form $\phi(t,R) = q(t)\f(R)$, where $q:[0,T)\to [0,\infty)$ and $\f:[0,1]\to [0,\infty)$ are both $C^2$.  We will also sometimes use the terminology \textit{homologous} synonymously with `separable.'

\begin{REMARK}
	Separable motions, by nature, allow for the separation of time and spatial dynamics, rendering them amenable to detailed study.  Their admissibility as solutions to \eqref{e:Eulerian} is tied to the symmetry of the system, or more precisely the constitutive assumptions that we make in Section \ref{ss:constitutive} below.  Separable motions also have a long history in the context of fluid and gas dynamics, which has influenced some of the more recent literature on elastic matter that is relevant to our work.  We will discuss this below, especially in Section \ref{ss:results}.
\end{REMARK}

\subsection{Basic Constitutive Assumptions}
\label{ss:constitutive}

Let $x(t,X)$ denote a flow map as defined in the previous subsection; let $Dx(t,X)$ denote its deformation tensor and $Jx(t,X) = \det Dx(t,X)$ its Jacobian determinant. The \textit{Piola-Kirchhoff stress tensor} $\sfS$ associated to the reference configuration $(\cB, \brho)$ is defined by
\begin{equation}
	\label{e:Sbardef}
	\sfS(t,X) = Jx(t,X) \sfT(t, x(t,X)) Dx(t,X)^{-T},
	\qquad (t,X)\in [0,T)\times \cB.
\end{equation}
Our constitutive assumptions will be stated in terms of $\sfS$ and related functions.  We need the following terminology, most of which can be found in standard references on mathematical elasticity theory, e.g., \cite{Gurtin, Ogden, MarsdenHughes}.
\begin{DEF}
	We say that the body $\cB(t)$ is \textit{elastic} if the Piola--Kirchhoff stress tensor $\sfS$ depends on time only through the deformation tensor, and if furthermore there exists a continuously differentiable `constitutive function' $\widehat{\sfS}:~GL_+(3,\R)\times~\cB\to~GL_+(3,\R)$ satisfying
	\begin{equation*} 
		\sfS(t,X) = \widehat{\sfS}(Dx(t,X),X), 
		\qquad X\in \cB, t\in [0,T).
	\end{equation*} 
	(Here $GL_+(3, \R)$ denotes the space of invertible $3\times 3$ matrices with real entries that have positive determinant.  We will also use the notation $SO(3,\R)$ to denote the subspace of $GL_+(3,\R)$ consisting of those matrices with determinant $1$.)	
	We say that $\cB(t)$ is \textit{hyperelastic} if it is elastic and, additionally, there exists a  \textit{strain-energy function} $W:GL_+(3,\R)\times \cB\to \R$ associated to $\widehat{\sfS}$, in the sense that 
	\begin{equation*}
		\widehat{\sfS}(F,X) = \brho \frac{\p}{\p F}W(F,X),
		\qquad F\in GL_+(3,\R), \;X\in \cB.
	\end{equation*} 
	We say that a hyperelastic body $\cB(t)$ is \textit{homogeneous} if the strain-energy function $W$ does not explicitly depend on $X$: 	
	\begin{equation*}
		W(F,X) = W(F).
	\end{equation*}	
\end{DEF}

\begin{DEF}
	We say that a strain-energy function $W$ is (i) objective or (ii) isotropic if
	\[
	\tag{i} W(F,X) = W(OF,X), \text{ for each } F\in GL_+(3,\R),\; O\in SO(3,\R),
	\]
	\[
	\tag{ii} W(F,X) = W(FO,X), \text{ for each } F\in GL_+(3,\R),\; O\in SO(3,\R).
	\]	
	We say that $W$ is homogeneous of degree $h$ if 	
	\[
	W(\s F) = \s^{h}W(F),
	\qquad \s>0, \;F\in GL_+(3,\R).
	\]
\end{DEF}

For the rest of the manuscript, we work with a hyperelastic, homogeneous body with a strain-energy function that is objective, isotropic, and homogeneous of degree $-1$.

\begin{REMARK}
	Hyperelasticity and homogeneity of the body, as well as objectivity and isotropy of $W$ are all mild assumptions.  Degree $-1$ homogeneity of $W$ is on the other hand a rather strong assumption.  However, as already remarked in Sideris's work \cite{Sideris2023} (Section 2), some homogeneity assumption is necessary in order to produce separable solutions---even in the absence of the gravitational term.  The gravitational force compels us to choose the degree of homogeneity to be equal to $-1$.  	
\end{REMARK}

\begin{REMARK}[Comparison with Mass-Critical Euler--Poisson] It is instructive to compare \eqref{e:Eulerian} with the Euler--Poisson system for a polytropic gas, particularly in the `mass-critical' case of polytropic exponent $\g = \frac43$.  This will be a frequent theme of the introduction.  
	
	If we set $W(F) = 3\k \big( \frac{\brho}{\det F} \big)^{\frac13}$ (which satisfies the constitutive hypotheses above for any $\k>0$), then the resulting Cauchy stress satisfies $\diver \sfT(t,x) = - \n(\k\rho(t,x)^{\frac43})$, in which case the system \eqref{e:Eulerian} becomes the Euler--Poisson equations for the case of a mass-critical polytropic gas, except for the boundary conditions, which are different than those typically imposed for Euler--Poisson system.  We refer to the boundary conditions \eqref{e:Eulerian}$_4$ as `solid vacuum boundary conditions,' to contrast with the `gaseous vacuum boundary conditions' where the pressure (and the density) is 
	required to vanish at the boundary.  
\end{REMARK}

In order to give $\diver \sfT$ the required units of force density, we assume that the constant $\k$ from the previous remark has the units of $\frac{\text{Force}}{\text{Length}^2 \cdot \text{Mass Density}^{\frac43}}$.  We normalize our (general, not-necessarily polytropic) $W$ and assume without loss of generality that 
\begin{equation} 
	\label{e:Wbarnormalization}
	W(I) = \brho^{\frac13}.
\end{equation}  
The apparent discrepancy in units between the stress term and the other terms in the equations of motion we will write down later can be resolved by taking into account the implicit units of $\k$.

\subsection{The Function $g$} 
\label{ss:g}
We now introduce a function that will be fundamental to the rest of the manuscript.  Under our working constitutive assumptions, we will show below in Section \ref{sss:W}
(c.f. also Lemma 6.1 of \cite{Sideris2023}) that the function 
\[
(\R_+, \mathbb{S}^{2}) \ni (y, \w) \mapsto W\big(I + (y-1)(\w \otimes \w)\big) 
\]
is independent of $\w$.  (Here $I$ denotes the $3\times 3$ identity matrix, and $\mathbb{S}^2$ is the $2$-sphere in $\R^3$.)  
This motivates us to define $g:(0,\infty)\to \mathbb{R}$ via 
\begin{equation} 
	\label{e:gdef} 
	g(y) = \brho^{- \frac13} W\big(I + (y-1)(\w \otimes \w)\big).
\end{equation} 

Sideris's analysis makes extensive use of a related quantity $f:(0,\infty)\to \R$, which (under the assumption of degree $-1$ homogeneity of $W$) is related to our $g$ via  $f(y) = y^{\frac13} g(y)$.  There is of course not much difference between $f$ and $g$ in terms of regularity, but working with $g$ rather than $f$ will allow us to write the equations of motion for a spherically symmetric flow map in a particularly compact form and will also simplify our derivation of this system.  

The constitutive hypotheses we have stated above will be sufficient to derive the equations of motion for spherically symmetric or 
separable solutions.  However, in order to actually prove the existence of our family of separable solutions, we will impose one additional constitutive assumption---on $g$ rather than on $W$ directly (c.f. \eqref{e:g''requirement} below).  Theorem 11.1 of \cite{Sideris2023} guarantees that this is mathematically legitimate.  More specifically, it states that if $f$ is strictly positive and sufficiently smooth, with $f'(1) = 0$ (equivalent to $g'(1) = -\frac13$), then there exists a strain-energy function $W$ that satisfies our other constitutive assumptions and such that (an identity that implies) \eqref{e:gdef} holds. 

\begin{REMARK}
	\label{rem:g(y)special}
	The case of a mass-critical polytropic gas corresponds (again, except for the boundary conditions) to the situation where $f(1) \equiv 1$, or equivalently, $g(y) = y^{-\frac13}$.  (This will become clear from the discussion in Section \ref{sss:W}.) 
\end{REMARK}

\subsection{Equations of Motion}
In Section \ref{s:Derivation}, we will derive equations for the flow map associated to a spherically symmetric solution of \eqref{e:Eulerian}, under the constitutive hypotheses discussed in Section \ref{ss:constitutive}.  As an easy consequence, we will obtain equations of motion for the separable solutions of interest, which we later prove exist (under one additional constitutive assumption).   We pause to record the relevant equations in the following Proposition and its Corollaries. 

\begin{PROP}
	\label{prop:spherical}
	Assume the constitutive hypotheses of the previous subsection hold (hyperelasticity and homogeneity of the body; objectivity, isotropy, and degree $-1$ homogeneity of the strain energy $W$).  Define $\cB$ as in \eqref{e:refdomain} and $g:(0,\infty)\to \R$ as in \eqref{e:gdef}.  Assume $\brho>0$.  Assume $\phi:[0,T)\times [0,1]\to \mathbb{R}$ is twice-differentiable in both arguments and that for all $t\in [0,T)$, $\phi$ satisfies $\phi(t,0)~=~\phi_{RR}(t,0)~=~0$, $\phi_R(t,\cdot)>0$.  Define (for all $t\in [0,T)$)
	\begin{equation}
		\label{e:Lambday}
		\Lambda(t,R) = \begin{cases} 
			\frac{\phi(t,R)}{R}, & R\in (0,1], \\
			\phi_R(t,0), & R = 0,
		\end{cases} 
		\qquad \qquad 
		y(t,R) = \frac{\phi_R(t,R)}{\Lambda(t,R)},
		\quad R\in [0,1].
	\end{equation}
	Assume that $\phi$ is a classical solution of the following system:
	\begin{equation}
		\label{e:spherical}
		\begin{cases}
			\displaystyle \phi_{tt} - \brho^{\frac13} \bigg[\frac{1}{R^2\phi} \p_R \bigg( \frac{R^4}{\phi} {g}'({y}) \bigg) + \frac{R}{\phi^2} {g}({y})\bigg] +  \frac{4\pi}{3} R^3 G \brho \phi^{-2} = 0,  \\
			{g}'({y}(t,1)) = 0.
		\end{cases}
	\end{equation} 
	Then there exists a corresponding classical solution $(\rho, \bu)$ of \eqref{e:Eulerian} on the time interval $[0,T)$.  The pair $(\rho, \bu)$ is related to the solution $\phi$ of \eqref{e:spherical} via $\rho(t, x(t, X)) \det Dx(t,X) = \brho$ and $\dot{x}(t,X) = \bu(t, x(t,X))$, where $x(t,X) = \phi(t,|X|)\frac{X}{|X|}$, $t\in [0,T)$, and $X\in \cB$.
\end{PROP}

In the special case where $\phi(t,R) = q(t)\f(R)$, it is straightforward to see that the equations of motion for $q$ and $\f$ can be decoupled, leading to the system of equations in the following Corollary.  

\begin{COROL} 
	\label{cor:separable}
	Assume the hypotheses of Proposition \ref{prop:spherical}. Assume additionally that $q:[0,T)\to [0,\infty)$ and $\f:~[0,1]\to~[0,\infty)$ are twice-differentiable, and that $\f$ satisfies $\f(0) = \f''(0) = 0$, $\f'>0$.  Define 
	\begin{equation} 
		\label{e:lambdaydef}
		\l(R) = \begin{cases} \frac{\f(R)}{R}, & R\in (0,1], \\
			\f'(0), & R = 0, \end{cases} 
		\qquad \qquad 
		y(R) = \frac{\f_R(R)}{\l(R)},
		\qquad R\in [0,1]. 
	\end{equation} 
	Assume that for some $\mu\in \R$, we have 
	\begin{equation}
		\label{e:separated}
		\begin{cases}
			\displaystyle 
			\frac{\f^2}{R^5} \p_R \bigg( \frac{R^4}{\f} g'(y)\bigg) + \frac{\f}{R^2}g(y) = \bigg( \frac{4\pi}{3}G \brho + \mu \l^3\bigg) \brho^{-\frac13} \f, \\
			g'(y(1)) = 0, \;\f(1) = 1, 
		\end{cases}
	\end{equation}
	and 
	\begin{equation}
		\label{e:q}
		q^2 \ddot{q} = \mu,
		\;\; q(0) = 1.
	\end{equation}	
	Then $\phi(t,R) = q(t)\f(R)$ satisfies \eqref{e:spherical}.	
\end{COROL} 

Finally, if we can guarantee that $g''(y)>0$, then we can rewrite \eqref{e:separated}$_1$ in the manner suggested by the following Corollary.
\begin{COROL}
	Assume the hypotheses of Corollary  \ref{cor:separable}.  If additionally, $g''(y(R))>0$ for all $R\in [0,1]$, $t\in [0,T)$, then \eqref{e:separated} is equivalent to the following:
	\begin{equation}
		\label{e:Reformulation1} 
		\begin{cases} 
			\frac{1}{R}\big( \f_{RR} + \frac{2}{R}( \f_R -\frac{\f}{R} ) \big) = - \frac{1}{R^2}( \f_R - \frac{\f}{R})\big[2(y-1) + \frac{E(y)}{g''(y)}\big] + \frac{\l}{\brho^{\frac13} g''(y)} \big( \frac{4\pi}{3} G \brho + \mu \lambda^3\big) \\
			g'(y(1)) = 0, \, \f(1) = 1.
		\end{cases} 
	\end{equation}
	where 
	\begin{equation}
		\label{e:defEh}
		E(y) = \begin{cases} 
			\frac{h(y) - h(1)}{y-1} - h'(y), & \; y \ne 1, \\
			0 & \;y = 1,
		\end{cases} 
		\qquad \qquad 
		h(y) = 3y g'(y) + g(y) = 3y^{\frac23} f'(y).
	\end{equation} 
\end{COROL}

\begin{REMARK}
	The additional assumption \eqref{e:g''requirement} that we make on $g$ below will guarantee not only that $g''(y)>0$ for the solutions we are interested in, but in fact that $g''(y)$ remains very large.  The special case $g(y) = y^{-\frac13}$ mentioned in Remark \ref{rem:g(y)special} does of course satisfy $g''(y)>0$, so that Corollary \ref{cor:separable} applies; however, $g(y) = y^{-\frac13}$ does not satisfy \eqref{e:g''requirement}, and thus our existence result (Theorem \ref{t:main} below) will not apply for this case.  
\end{REMARK}

\begin{REMARK}
	\label{rem:q}
	In what follows, our main focus will be Equations \eqref{e:separated} and \eqref{e:Reformulation1}, which are nontrivial to analyze.  On the other hand, Equation \eqref{e:q} is elementary \cite{FuLin1998, HadzicJang2018CPAM}; its behavior depends on the sign of $e_{\eff}~=~\frac12 \dot{q}(0)^2 + \mu$.  If $e_\eff>0$, then $q(t)$ grows at an asymptotically linear rate for large $t$; we refer to this as the `linear expanding' case.  
	If $e_{\eff} = 0$ (which can only occur when $\mu\le 0$), we have $\dot{q}(0) = \pm \sqrt{-2\mu}$ and $q(t) = (1 + \frac32 \dot{q}(0)t)^{\frac23}$ for all $t\in [0,T)$, which can lead to  $q\equiv 1$ (the `stationary' case) or expansion or collapse, depending on the sign of $\dot{q}(0)$.
	The case $e_{\eff}<0$ (which can only occur when $\mu<0$) always corresponds to finite-time collapse, with $q(t){\sim}_{t\to T^-} \;c_1(T - c_2 t)^{\frac23}$.
\end{REMARK}

\subsection{Main Result and Relation to Existing Literature}

\label{ss:results}

We now state our main result.
\begin{THEOREM}
	\label{t:main}
	Assume the hypotheses of Proposition \ref{prop:spherical} hold.  Assume furthermore that 	$g:~(0,\infty)\to~\mathbb{R}$ is $C^3$ and that 
	\begin{equation}
		\label{e:g''requirement}
		g''(1) \ge 50 \bigg( 50 + \sup_{|y - 1|\le \frac12} |g'''(y)|\bigg).
	\end{equation}
	Then there exists $\mu_0>0$ such that for every $\mu\in [-\mu_0, \mu_0]$, there exists $\brho>0$ and a corresponding classical solution $\f$ 
	to the system \eqref{e:separated} equipped with the parameters $\mu, \brho$.
\end{THEOREM}

Theorem \ref{t:main} provides the existence of homologous expanding solutions, collapsing solutions, and stationary solutions for self-gravitating hyperelastic matter with degree $-1$ homogeneity.  Whether a solution experiences expansion or collapse, or remains stationary, depends on the values of $\mu$ and $e_{\mathrm{eff}}$ (c.f. Remark \ref{rem:q}). To the best of our knowledge, the finite-time collapsing solutions provided by Theorem~\ref{t:main} constitute the first such rigorous result in the presence of positive residual pressure\footnote{Following \cite{Sideris2023}, we define the \textit{residual pressure} $\cP(1)$ via $\cP(1)I = -\widehat{S}(I)$.  It follows from the discussion in \ref{sss:g1g'1} below that $\cP(1) = \frac13 \brho^{\frac43}$ in our setting.}.  We will prove Theorem \ref{t:main} in Section \ref{s:proof}, after providing a concise derivation of \eqref{e:spherical} from \eqref{e:Eulerian} in Section \ref{s:Derivation}. 

The most closely related results in the literature to our Theorem \ref{t:main} appear in two recent papers by  Sideris \cite{Sideris2023} and Calogero \cite{Calogero2022}.  Let us briefly summarize the relevant findings of those two works as they relate to ours. 
Sideris \cite{Sideris2023} considers separable dynamics \textit{in the absence of self-gravitation}.  Whether Sideris's solutions expand or collapse depends entirely on the sign of the residual pressure; positive residual pressure corresponds to expansion, while negative residual pressure corresponds to finite-time collapse.  Adding gravity to the system in our present work provides a mechanism to compete with the positive residual pressure and allows for the possibility of collapsing solutions in a qualitatively new parameter range.  The reader familiar with \cite{Sideris2023} will notice a great deal of similarity between the framework outlined in Sideris's paper and the one we use in our work, and our analysis indeed owes a great debt to Sideris's. Nonetheless, 
we still feel it is valuable to develop our framework from the ground up, because our use of the function $g$ rather than Sideris's $f$ substantially simplifies equations \eqref{e:spherical} and \eqref{e:separated} and streamlines their derivations (as well as a good deal of the subsequent analysis).

In \cite{Calogero2022}, Calogero studies what he refers to as the `polytropic elastic model.'  The main contribution of that paper is to extend the constitutive theory for polytropic gases in a natural way to apply to spherically symmetric elastic bodies.  Notably, \cite{Calogero2022} is set entirely in the Eulerian formulation.  As an application of his analysis, Calogero considers homologous dynamics similar to those treated in our paper.  Some of his constitutive assumptions are quite different from ours, and the results of \cite{Calogero2022} do not appear to be directly comparable to those of the present work.  However, \cite{Calogero2022} does establish the existence of a certain family of self-gravitating separable expanding solutions with solid vacuum boundary conditions, and it also provides numerical evidence for the existence of separable collapsing solutions with gaseous vacuum boundary conditions. We also mention the works \cite{AC2020, ACL2022}, where the authors studied self-gravitating spherically
symmetric elastic bodies in static equilibrium. 

\begin{REMARK}
	The assumption \eqref{e:g''requirement} is a stringent requirement and very far from optimal, but some largeness of $g''(1)$ does appear to be needed, and the particular choice \eqref{e:g''requirement} helps to streamline parts of the argument below.  Note also that for convenience, we will routinely write down non-sharp estimates in the following analysis as well.
\end{REMARK}

\begin{REMARK}
	The condition \eqref{e:g''requirement} easily implies the following one, which is written in a way that reflects how we will actually use the largeness of $g''(1)$ below.  
	\begin{equation}
		\label{e:g''requirement.mod}
		g''(1) \ge 2400(1 + M^2), 
		\qquad \qquad M = \sup_{|y-1|\le \frac12} \frac{|g'''(y)|}{g''(1)}.
	\end{equation}
\end{REMARK}

We close our introduction by comparing our work with some relevant literature on fluid and gas dynamics.  For a self-gravitating gas with mass-critical polytropic exponent $\g = \frac43$ (i.e., the special case $g(y) = y^{-\frac13}$, with gaseous vacuum boundary conditions), the existence of homologous solutions goes back to the works of Goldreich--Weber \cite{GoWe1980}, Makino~\cite{Makino1992}, and Fu--Lin \cite{FuLin1998}.  The affine ansatz separates the time and spatial dynamics: the time dynamics are described by precisely the same ODE \eqref{e:q} as in our setup, while the background deformation profile is determined by the generalized Lane--Emden equation, written in terms of the \textit{enthalpy} profile $w = (\rho^0)^{\g-1} =~(\rho^0)^{\frac13}$:
\begin{equation}
	\label{e:genLE}
	\begin{cases} 
		w'' + \frac{2}{r} w' + \pi G w^3 = -\frac{3}{4} \mu   \\
		w'(0) = 0, \;\; w(1) = 0.
	\end{cases}
\end{equation}
(The boundary conditions \eqref{e:genLE}$_2$ correspond to the formulation in \cite{HadzicJang2018CPAM}, which is easiest to compare to our present setup.) An essentially complete analysis on the enthalpy distribution satisfying this semi-linear elliptic equation and gaseous vacuum boundary condition is available for the full allowed range of $\mu$.  In particular, there exists a maximal curve of solutions to \eqref{e:genLE} in the $(\mu, \brho)$ parameter space, with $\brho\in (0, \brho_*)$; furthermore, there exists a critical value $\brho_c\in (0,\brho_*)$ below which only expanding solutions exist, and above which collapsing solutions are admitted. See \cite{FuLin1998} for more details.  The expanding solutions of \cite{GoWe1980, Makino1992, FuLin1998} were shown to be nonlinearly stable in \cite{HadzicJang2018CPAM}, indicating that this type of expansion is a stable phenomenon and provides a mechanism for avoiding compressive singularities. See also \cite{Liu2020} for the stability of expanding solutions for viscous radiation gaseous stars.  In forthcoming work, the authors plan to prove that the expanding solutions from Theorem \ref{t:main} are also nonlinearly stable.

The analog of \eqref{e:genLE} for a self-gravitating hyperelastic material is as follows: 
\begin{equation}
	\label{e:genLEanalog}
	\begin{cases} 
		-\frac{3}{4r^2} \partial_r \bigg( \frac{1}{w^3 r} \p_r \big[ w^4 r^3 z^{\frac43} g'(z) \big] + r w z^{\frac13} g(z) \bigg) + \pi G w^3 = -\frac{3}{4}\mu \\
		z(r) = \frac{3}{w^3} \int_0^r w(s)^3 s^2 \ds  \\
		w'(0) = 0, \;\; g'(z(1)) = 0.
	\end{cases}
\end{equation}
Note that Equation \eqref{e:genLE}$_1$ can be recovered from \eqref{e:genLEanalog}$_1$ by taking $g(y) = y^{-\frac13}$, in which case $z$ drops out.  Note also the difference in the boundary conditions between the two systems.  It would be 
interesting to carry out an investigation for self-gravitating hyperelastic matter analogous to what has been done for \eqref{e:genLE} in the context of self-gravitating gas, for instance as in \cite{FuLin1998}.  

Our analysis below takes a very different approach,  
which is essentially a perturbative argument: it is based on  the fixed point theorem for the integral reformulation of \eqref{e:Reformulation1} (cf. \eqref{e:FP_Reformulation}--\eqref{e:cFdef}) for suitable range of $(\mu, \brho)$ and a continuity argument for matching the boundary condition.   
We make no attempt to obtain a \textit{maximal} curve of solutions in $(\mu, \brho)$ parameter space.  

As a final note, we mention that in the fluid and gas literature (see for example \cite{LiWang2006, Ovs1956}), the admissible separable dynamics include spherically symmetric expansion/collapse as well as more general affine motions. In the context of isolated moving gases, Sideris in \cite{Sideris2017Affine} constructed affine solutions to the vacuum free boundary Euler equations representing moving ellipsoids that expand and spread to infinity with a linear rate, and moreover in \cite{HadzicJang2018Inventiones, ShSi2019}, Sideris’s solutions were shown to be nonlinearly stable. 

\section{Derivation of the Equations of Motion for Separable Solutions}
\label{s:Derivation}

\subsection{General Lagrangian Formulation}

Given a solution $(\rho, \bu)$, we define the flow map $x$ as before, namely a function $x:[0,T)\times \cB \to \cB(t)$ that satisfies the relations $\dot{x}(t,X) = \bu(t,x(t,X))$ and $\rho(t, x(t,X)) \det Dx(t,X) = \brho$, for all $t\in [0,T)$ and all $X\in \cB$. 
We set the following notation:
\[
\begin{split} 
	\Psi(t,X) & = \Phi[\rho(t)](t,x(t,X)), \\
	A(t,X) & = Dx(t,X)^{-1}. 
\end{split} 
\]

We use the notation $a_{,k}$ to denote the derivative of a scalar quantity $a(t,\cdot)$ with respect to its $k$th argument, excluding time.

Setting $x = x(t,X)$ in \eqref{e:Eulerian}$_2$ and multiplying by $Jx(t,X)$, we obtain 
\begin{equation} 
	\label{e:Langrangian1bar}
	\brho \ddot{x}^i = \sfS^i_{j,j} - G\brho A_i^j\Psi_{,j},
	\qquad \text{ in } \cB.
\end{equation} 
Here we have used 
\[
\sfS^i_{j,j}(t,X) = Jx(t,X) \sfT^i_{j,j}(t, x(t,X)),
\]
which follows from the definition \eqref{e:Sbardef} of $\sfS$ and the Piola identity
\begin{equation} 
	\label{e:Piolaid}
	(A^j_k Jx )_{,j} = 0.
\end{equation} 
The full system in these Lagrangian coordinates thus becomes 
\begin{equation} 
	\label{e:Langrangian}
	\begin{cases} 
		\brho \ddot{x}^i = \sfS^i_{j,j} - G\brho A_i^j\Psi_{,j},
		& \qquad \text{ in } [0,T)\times \cB. \\
		A^k_i(A^\ell_i \,\Psi_{,\ell})_{,k} = 4\pi \brho (Jx)^{-1}, & \qquad \text{ in } [0,T)\times \cB, \\
		\sfS(t,X)\bn(X) = 0, & \qquad X\in \p \cB.
	\end{cases}
\end{equation}	

\subsection{Lagrangian Formulation for Spherically Symmetric Solutions}

Under the assumption of spherical symmetry (and using our constitutive hypotheses), we can simplify \eqref{e:Langrangian} substantially.  
We write $x(t,0)~=~0$ and 
\begin{equation}
	\label{e:xphi}
	x(t,X) = \phi(t,R)\w,
	\qquad R = |X|,\;\;\w = \frac{X}{R},
	\qquad X\ne 0.
\end{equation} 
We also assume (as above) that 
\[
\phi(t,0) = \phi_{RR}(t,0) = 0, \text{ for } t\ge 0;
\qquad 	\phi_R(t,R)>0 \text{ for } t\ge 0,\, R\in [0,1],
\]
which guarantees that $x(t, \cdot)$ is $C^2$ even at the origin.  

\subsubsection{Deformation Tensor and Strain-Energy Function}
\label{sss:W}
A short computation shows that the deformation tensor takes the form 
\begin{equation}
	\label{e:deftensor1} 
	Dx(t,X) = \Lambda(t,R) \big( I + (y(t,R) - 1)(\w \otimes \w) \big), 
\end{equation}
where $\Lambda(t,R)$ and $y(t,R)$ are defined as in \eqref{e:Lambday}.  
Note that the eigenvalues of $Dx(t,X)$ are $\phi_R(t,R)$ (with 1-dimensional eigenspace $\Span\{\w\}$) and $\L(t,R)$ (with 2-dimensional eigenspace $(\Span\{\w\})^\perp$).  Motivated by \eqref{e:deftensor1}, we introduce the notation 
\[
F(y,\L, \w) = \L \big( I + (y-1)(\w\otimes \w)\big),
\qquad y,\L>0,\;\w\in \S^2.
\]
Note that if $O\in SO(3,\R)$, then clearly
\[
O F(y, \L, \w) O^T  = F(y,\L, O\w).
\]
Therefore, objectivity and isotropy of $W$ guarantee that 
\begin{equation}
	W(F(y,\L, \w)) = W(F(y,\L, O\w)),
	\qquad O\in SO(3,\R).
\end{equation}
It follows that $W(F(y,\L, \w))$ is independent of $\w = \frac{X}{R}$, so we may define 
\[
g(y) = \brho^{-\frac13} W\big( F(y,1,\w) \big),
\qquad y>0,
\]
as claimed in Section \ref{ss:g}.

\subsubsection{The Piola--Kirchoff Stress Tensor}

As pointed out by Sideris \cite{Sideris2023}, it is well-known (c.f. \cite{Ogden}, Theorem 4.2.5) that for a hyperelastic, homogeneous material with an objective, isotropic strain-energy function $W$, one has
\[
\sfT(t,x(t,X)) \in \Span\{I, \;(Dx(t,X) Dx(t,X)^T)^{\frac12}, \; Dx(t,X) Dx(t,X)^T \}.
\]
Under our assumptions of spherical symmetry, together with the relationship between $\sfS$ and $\sfT$, it follows that 
\[
\widehat{\sfS}(Dx(t,X)) \in \Span\{ P(\w), Q(\w)\},
\] 
where $P(\w) = \w\otimes \w$ and $Q(\w) = I - P(\w)$.  We may therefore write 
\begin{equation}
	\label{e:Sbarspan}
	\widehat{\sfS}(F(y,\L,\w)) = c_1(y,\L, \w) P(\w) + c_2(y,\L, \w) Q(\w),
\end{equation}
where $c_1$ and $c_2$ are determined as follows.  Noting that $P^i_j P^i_j = 1$, $P^i_j Q^i_j = 0$, and $Q^i_j Q^i_j = 2$, we have
\begin{align*}
	\brho^{\frac43} \L^{-1} {g}'(y)
	& = \brho \frac{\p}{\p y} {W}(F(y,\L,\w)) = {\sfS}^i_j \L P^i_j = \L c_1, \\ 
	-\brho^{\frac43} \L^{-2} {g}(y)
	& = \brho \frac{\p}{\p \L} {W}(F(y,\L, \w)) =  {\sfS}^i_j (yP^i_j + Q^i_j) = y c_1 + 2c_2.
\end{align*}
Thus 
\begin{equation}
	c_1(y,\L,\w) = \brho^{\frac43}\L^{-2} {g}'(y),
	\qquad 
	-2c_2(y,\L,\w) = \brho^{\frac43}\L^{-2}\big( y {g}'(y) + {g}(y) \big).
\end{equation}
Notice that both $c_1$ and $c_2$ are independent of $\w$. We record the formula for $\widehat{\sfS}$ for further use below:
\begin{equation}
	\label{e:Shatformula}
	\brho^{-\frac{4}{3}} \widehat{\sfS}(F(y,\L,\w)) = \L^{-2} {g}'(y) P(\w) -\frac12 \L^{-2}\big( y {g}'(y) + {g}(y) \big) Q(\w).
\end{equation}

\subsubsection{Divergence of the Piola--Kirchhoff Stress}

We now substitute $y = {y}(t,R)$, $\L = {\L}(t,R)$ into \eqref{e:Sbarspan} and take the divergence (in $X$), recalling that $c_1$ and $c_2$ are independent of $\w$. We obtain 
\begin{align*}
	\w_i \sfS^i_{j,j}(t,X) & = \p_R(c_1({y},{\L})) + \frac{2}{R}(c_1({y},{\L}) - c_2({y},{\L})) 
	= \frac{1}{R^2} \p_R \big( R^2 c_1 ({y},{\L}) \big) - \frac{1}{R} \cdot 2c_2({y},{\L}).
\end{align*}

Substituting in the formulas for $c_1$ and $c_2$, we obtain 
\begin{align*}
	\brho^{-\frac43}\w_i \sfS^i_{j,j}(t,X) 
	& = \frac{1}{R^2} \p_R \bigg( \frac{R^2}{\L^2} {g}'({y}) \bigg) + \frac{1}{R \L^2} ({y} {g}'({y}) + {g}({y})) 
	= \frac{1}{R^3 \L} \p_R \bigg( \frac{R^3}{\L} {g}'({y}) \bigg) + \frac{1}{R \L^2} {g}({y}).
\end{align*} 

Finally, substituting $\L(t,R) = \frac{\phi(t,R)}{R}$, we obtain the following expression (with the understanding that $\phi(R)/R$ should be replaced by $\phi'(0)$ when $R = 0$): 
\begin{equation}
	\label{e:divSfinal}
	\brho^{-\frac43} \w_i {\sfS}^i_{j,j}(t,X) 
	= \frac{1}{R^2\phi} \p_R \bigg( \frac{R^4}{\phi} {g}'({y}) \bigg) + \frac{R}{\phi^2} {g}({y})
\end{equation}

\subsubsection{The Potential Term}

Following (for example) \cite{Jang2014}, the Poisson equation yields that
\begin{equation}
	\label{e:Potential}  
	\w^i A^i_j \Psi_{,j} = \frac{1}{\phi^2} \int_0^R 4\pi \brho s^2 \ds
	= \frac{4\pi}{3} R^3 \brho \phi^{-2}.
\end{equation} 

\subsubsection{Boundary Conditions}
Our boundary conditions read ${\sfS}(t,X)\w = 0$ for $X\in \p \cB$.  Equivalently, since $P(\w)\w = \w$ and $Q(\w)\w = 0$, our boundary conditions can be rewritten as
\begin{align*}
	0 
	& = \widehat{\sfS}(F({y}(t,1), \L(t,1), \w))\w = \brho^{\frac43} \L(t,1)^{-2} g'(y(t,1)),
\end{align*}
which reduces to \eqref{e:spherical}$_2$.

\subsubsection{Full Lagrangian System}

Multiplying \eqref{e:Langrangian}$_1$ by $\w_i$, then substituting \eqref{e:xphi}, \eqref{e:divSfinal}, and \eqref{e:Potential} yields exactly \eqref{e:spherical}$_1$.  This completes the derivation of the system \eqref{e:spherical}.

\subsubsection{Determining $g(1)$ and $g'(1)$} 

\label{sss:g1g'1}

The definition of $g$, together with our normalization \eqref{e:Wbarnormalization}, imply that 
\begin{equation} 
	\label{e:g(1)}
	g(1) = \brho^{-\frac13} W(F(1,1,\w)) = \brho^{-\frac13} W(I) = 1.
\end{equation} 
To determine $g'(1)$, we note that substituting $y = 1$ into \eqref{e:Shatformula} yields (after a short computation, and substituting $g(1) = 1$)
\[
\brho^{-\frac43} \widehat{\sfS}(\Lambda I) = \frac{1}{2\Lambda^2} \big( (1 + 3g'(1))(\w \otimes \w) - (1 + g'(1))I \big).
\]
Since $\widehat{\sfS}$ is required to be continuously differentiable in $X$ by assumption (in particular, at $X=0$), the first term on the right side of the above expression must vanish, forcing 
\begin{equation} 
	\label{e:g'(1)}
	g'(1) = -\frac13.
\end{equation}

\subsubsection{Derivation of \eqref{e:Reformulation1}}

We now work with a separable solution and write $\phi(t,R) = q(t)\f(R)$.  We will also use the notation for $\l(R)$ and $y(R)$ introduced in \eqref{e:lambdaydef}.  Deriving \eqref{e:separated} and \eqref{e:q} from \eqref{e:spherical} is straightforward, and we leave the details to the reader.  We do, however, fill in the remaining details necessary to get from \eqref{e:separated} to \eqref{e:Reformulation1}, under the additional assumption that $g''(y(R))>0$ for all $R\in [0,1]$.  We rewrite the definitions of $E$ and $h$ below for the reader's convenience:
\[
E(y) = \begin{cases} 
	\frac{h(y) - h(1)}{y-1} - h'(y), & \; y \ne 1, \\
	0 & \;y = 1,
\end{cases} 
\qquad \qquad 
h(y) = 3y g'(y) + g(y).
\]
Note that $h(1) = 0$ by \eqref{e:g(1)} and \eqref{e:g'(1)}. 
Note also that if $g$ is $C^k$, then $h$ is $C^{k-1}$ and $E$ is $C^{k-2}$.  


We now rewrite the left side of \eqref{e:separated}$_1$ as follows.  
\begin{align*}
	\frac{\f^2}{R^5} \p_R \bigg( \frac{R^4}{\f} g'(y)\bigg) + \frac{\f}{R^2}g(y) 
	& =  g''(y) \f_{RR} + \frac{\f}{R^2}  \big[ y (1-y)g''(y) + (4-y)g'(y)  + g(y) \big] \\
	& =  g''(y) \f_{RR} + \frac{\f}{R^2}  \big[ -2y(1-y)g''(y) + (1-y) h'(y) +  h(y) \big] \\
	& = g''(y) \bigg( \f_{RR} + \frac{2\f}{R^2}(y-1) \bigg)  
	+ \frac{\f}{R^2} (y-1) \big[ 2(y-1) g''(y) + E(y) \big] 
\end{align*}

Under our assumption that $g''(y)>0$, we may divide by $R g''(y)$ and rewrite the whole of equation \eqref{e:separated}$_1$ as exactly \eqref{e:Reformulation1}$_1$.

\section{Proof of Theorem \ref{t:main}} 

\label{s:proof}

The rest of the paper is dedicated to the proof of our main result, Theorem \ref{t:main}.  We provide a brief outline of our proof: In Section \ref{ss:FP}, we write the eigenvalue problem \eqref{e:Reformulation1}$_1$ as a fixed point problem of a certain operator, parametrized by the putative eigenvalue $\mu$ and reference density $\brho$.  In Section \ref{ss:Domaing}, we restrict the domain of $g$ to guarantee that $g''(y)>0$ (so that the use of \eqref{e:Reformulation1} is justified) and provide some preliminary estimates.  In Section~\ref{ss:Contraction}, we give conditions under which the operator introduced in Section~\ref{ss:FP} is a contraction mapping and thus has a unique fixed point.  In Section \ref{ss:BC}, we fix $\mu$ and show that we can tune $\brho$ to obtain a fixed point of our operator that \textit{also satisfies the boundary conditions}; the latter are not fully taken into account in the previous steps.

\subsection{Fixed Point Problem}
\label{ss:FP}
Our starting point is Equation \eqref{e:Reformulation1}, which we rewrite as 
\begin{align}
	\label{e:Reformulation2} 
	\frac{1}{R}\bigg( \f_{RR} + \frac{2}{R}\big( \f_R -\frac{\f}{R} \big) \bigg) 
	& = - \frac{1}{R^2}\big( \f_R - \frac{\f}{R}\big)U(y) + \frac{V(\l)}{g''(y)},
\end{align}
where 
\begin{equation} 
	\label{e:UVdef}
	U(y) = 2(y-1) + \frac{E(y)}{g''(y)},
	\qquad 
	V(\lambda) = \frac{\l}{\brho^{\frac13}} \bigg( \frac{4\pi}{3} G \brho + \mu \lambda^3\bigg).
\end{equation} 

Note that $U(1) = 0$, and that $U$ and $V$ are differentiable.  We will provide estimates on their derivatives later on.  First, we give one last reformulation of \eqref{e:Reformulation2}.

In order to enforce $\f(0) = \f''(0) = 0$, we make the ansatz $\f''(R) = R \zeta(R)$, where $\zeta\in C([0,1])$.  By a short computation, we thus have
\begin{equation}
	\label{e:phifromzeta}
	\f(R) = R\f'(0) + \int_0^R (R - \t)\t \zeta(\t)\dtau.
\end{equation} 
The normalization $\f(1) = 1$ forces
\[
\f'(0) = 1 - \int_0^1 (1 - \t) \t \zeta(\t) \dtau.
\]
However, we will continue to write $\f'(0)$ for brevity.  We also have 
\[
\f(R) \equiv R 
\qquad 
\text{ when } 
\qquad \zeta \equiv 0.
\]
The following are also useful:
\begin{align*} 
	\f'(R) & = \f'(0) + \int_0^R \t \zeta(\t)\dtau \\
	\f'(R) - \lambda(R) & = \frac1R \int_0^R \t^2 \zeta(\t) \dtau  \\
	\frac1R \bigg( 
	\f_{RR} + \frac{2}{R}\big( \f_R - \lambda \big) \bigg) 
	& = \zeta(R) + \frac{2}{R^3} \int_0^R \t^2 \zeta(\t) \dtau.
\end{align*} 

We may thus rewrite our entire system as 
\begin{equation} 
	\label{e:FP_Reformulation}
	L\zeta = \cF[\zeta],
\end{equation} 
where 
\begin{equation} 
	\label{e:Lzetadef}
	L\zeta(R) = \zeta(R) + \frac{2}{R^3} \int_0^R \t^2 \zeta(\t) \dtau,
\end{equation} 
\begin{align}
	\label{e:cFdef}
	\cF[\zeta](R) 
	& = -\frac{1}{R^2} \big( \f_R - \l \big) U(y) + \frac{V(\l)}{g''(y)}.
\end{align}
Note that $\f$ and $\f'$ (and thus $\lambda$ and $y$) are determined by $\zeta$, so the right side of \eqref{e:cFdef} is determined by $\zeta$ as well. Note furthermore that $L:C([0,1])\to C([0,1])$ is invertible, and we can find its inverse explicitly.  (Write $L\zeta = \eta$, multiply by $R^3$, differentiate, and solve for $\zeta$.) The formula is:
\[
L^{-1} \eta(R) = \eta(R) - \frac{2}{R^5} \int_0^R \t^4 \eta(\t)\dtau.
\]
Our task of finding a solution $\f$ of \eqref{e:separated}$_1$ is now reduced to that of finding a fixed point $\zeta$ of $\zeta \mapsto L^{-1}\cF[\zeta]$.  For use below, we record that 
\begin{equation}
	L^{-1}[1] = \frac35, 
	\qquad \|L^{-1}\|_{op} = \frac75.
\end{equation} 
(It is clear that $\|L^{-1}\|_{op}\le \frac75$; to see the opposite inequality, consider the sequence of functions defined by $\eta_n(R) = -1 + 2R^n$.)

\subsection{Restricted Domain of $g$} 
\label{ss:Domaing}
The equations \eqref{e:FP_Reformulation}, \eqref{e:Lzetadef}, \eqref{e:cFdef} constitute a valid reformulation of \eqref{e:separated}$_1$ if we can guarantee that $g''(y(R)) > 0$ for all $R\in [0,1]$.  We now restrict the domain of $g$ to an interval $I(\d)$ on which $g''$ is positive, justifying later that $y(R)$ actually belongs to this interval for the solutions we consider.  Under the assumption \eqref{e:g''requirement}, we define 
\begin{equation} 
	\label{e:delta} 
	\delta = \frac{10}{g''(1)}.  
\end{equation} 
Then 
\begin{equation}
	\frac{9}{10} g''(1) \le g''(y)\le \frac{11}{10} g''(1)
	\qquad 
	\text{ whenever } 
	\qquad |y - 1| \le \d.
\end{equation}
We will work on the fixed interval $I(\delta) = [1 - \delta, 1+\delta]$.  It is not difficult to obtain the following estimates.  (Here $M$ is as defined in \eqref{e:g''requirement.mod}.)
\[
\|h''\|_{L^\infty(I(\d))} \le 7\big( \|g'''\|_{L^\infty(I(\d))} + g''(1) \big), 
\qquad 
\|E'\|_{L^\infty(I(\d))} \le 2 \|h''\|_{L^\infty(I(\d))}.
\]

\begin{equation}
	\label{e:Uinfo} 
	\|U'\|_{L^\infty(I(\d))} \le 2 + \frac{2\|h''\|_{L^\infty(I(\d))}}{\frac{9}{10}g''(1)} + \frac{2 \|h''\|_{L^\infty(I(\d))} |y-1| M}{\frac{9}{10} g''(1)} \le 22 ( 1 + M^2).
\end{equation} 
Note that we also have
\begin{equation}
	\label{e:Vinfo} 
	|V(\l)|\le (1 + |\l - 1|)^4 K(\brho, \mu),
	\qquad 
	|V'(\l)|\le 4(1 + |\l - 1|)^3 K(\brho, \mu),
\end{equation} 
where 
\begin{equation} 
	\label{e:Kbrhomu}
	K(\brho, \mu) := \brho^{-\frac13} \big( \frac{4\pi}{3} G \brho + |\mu|\big).
\end{equation}

\subsection{The Contraction Mapping}
\label{ss:Contraction}

We show that, under additional assumptions on $\brho$ and $\mu$, the operator $L^{-1} \cF$ is a contraction on $N(\d)$, where 
\[
N(\d):=\{ \zeta\in C([0,1]): \|\zeta\|_\infty\le \d\},
\] 
and where $\d>0$ is defined as in \eqref{e:delta}. 

\subsubsection{Basic Estimates}

Given $\zeta\in N(\d)$, let $\f$ be defined via \eqref{e:phifromzeta}, and let $\l$, $y$ be defined as in \eqref{e:lambdaydef}.  Let $\widetilde{\zeta}$ be another element of $N(\d)$ and define $\widetilde{\f}$, $\widetilde{\l}$, $\widetilde{y}$ in the analogous way.  Then after a few short computations, we may estimate (for $R\in [0,1]$):
\begin{equation} 
	|\f_R(R) - \widetilde{\f}_R(R)| \le \frac{R^2}{2}\|\zeta - \widetilde{\zeta}\|_\infty,
	\qquad 
	|\l(R) - \widetilde{\l}(R)| 
	\le \frac{R^2}{6} \|\zeta - \widetilde{\zeta}\|_{\infty},\qquad 
	\zeta, \widetilde{\zeta} \in N(\d).
\end{equation} 
\begin{equation} 
	|y(R) - \widetilde{y}(R)|
	\le R^2\|\zeta - \widetilde{\zeta}\|_{\infty},  
	\qquad 
	\zeta, \widetilde{\zeta} \in N(\d).
\end{equation} 
Of course, if $\widetilde{\zeta} \equiv 0$, these reduce to
\begin{equation}
	|\f_R(R) - 1| \le \frac{R^2 \d}{2},
	\qquad 
	|\l(R) - 1|
	\le \frac{R^2 \d}{6},
	\qquad 
	|y(R) - 1|
	\le R^2 \d, \qquad \qquad\qquad  \zeta \in N(\d).
\end{equation} 
The last estimate on the right confirms that  $y(R)\in I(\d)$ whenever $\zeta\in N(\d)$.

\subsubsection{Range of $\mu$} In order to obtain proper estimates on $\cF$, we restrict the range of $\brho$ and $\mu$.

First, fix $\brho_+>0$ satisfying $\frac{4\pi}{3} G \brho_+^{\frac23} = 10$.  In other words, put
\begin{equation} 
	\label{e:brho+}
	\brho_+ = \bigg( \frac{10}{\frac{4\pi}{3} G} \bigg)^{\frac32}.
\end{equation} 

Next, note that for fixed $\mu$, the function $\bar{\varrho} \mapsto K(\bar{\varrho}, \mu)$ is minimized when $\bar{\varrho} = \brho_-(\mu)$, where 
\begin{equation} 
	\label{e:brho-}
	\brho_-(\mu) := \frac{|\mu|}{\frac{8\pi}{3}G}.
\end{equation} 
(Note: When $\mu = 0$, \eqref{e:brho-} implies that $\brho_-(\mu) = 0$.  We stress that $\bar{\varrho} = 0$ is a valid input for the function $K(\cdot, \mu)$, even as we continue to insist that the parameter $\brho$ that appears in our Lagrangian formulation is always strictly positive.) Furthermore, the minimum value $K(\brho_-(\mu), \mu)$ is equal to some absolute constant times $|\mu|^{\frac23} G^{\frac13}$.  

With the above in mind, we may pick $\mu_0>0$ small enough so that both the following conditions are satisfied:
\begin{itemize}
	\item First, we require that $\mu_0 < \min\{\frac12\brho_+^{\frac13}, \frac{8\pi}{3} G \brho_+\}$.  This has the following three implications: 
	\begin{itemize}
		\item The fact that $\mu_0 < \frac{8\pi}{3} G \brho_+$ implies that 
		\begin{equation} 
			\label{e:brho-<brho+}
			\brho_-(\mu) < \brho_+,
			\qquad \qquad \text{ whenever } |\mu|\le \mu_0. 
		\end{equation} 
		Thus $\brho\mapsto K(\brho, \mu)$ is increasing in $\brho$ for $\brho\in [\brho_-(\mu), \brho_+]$.
		\item The fact that $\mu_0 < \frac12 \brho_+^{\frac13}$ implies that
		\begin{equation}
			0 \le K(\brho, \mu) \le K(\brho_+, \mu) < \frac{21}{2}
			\qquad \text{ whenever } |\mu|\le \mu_0 \text{ and } \brho\in [\brho_-(\mu), \brho_+]. 
		\end{equation}
		The fact that $\mu_0 < \frac12\brho_+^{\frac13}$ also implies that 
		\begin{equation} 
			\label{e:epsilonlwr}
			\epsilon(\brho_+, \mu) := \frac{4\pi}{3} G \brho_+^{\frac23} + \mu \brho_+^{-\frac13}> \frac{19}{2},
			\qquad \text{ whenever } |\mu|\le \mu_0.
		\end{equation} 
	\end{itemize}
	The first two of these implications will be put to almost immediate use; the utility of \eqref{e:epsilonlwr} will become apparent when we discuss boundary conditions.  
	
	\item Second, we require that 
	\begin{equation} 
		\label{e:Ksmall}
		K(\brho_-(\mu), \mu) < \frac{1}{20},
		\qquad 
		\text{ whenever } |\mu|\le \mu_0.
	\end{equation} 
	This requirement will also be used only when we discuss boundary conditions.
\end{itemize}

\begin{REMARK}
	\label{r:nonoptimality}
	The reader may observe that it is possible to expand the range of $\mu$ beyond the interval $[-\mu_0, \mu_0]$ that we cover.  For example, more precise bounds on $V$ and $V'$ than \eqref{e:Vinfo} are available for negative values of $\mu$, and by using such bounds and modifying the argument below, existence of separable solutions can be established for certain values of $\mu$ outside $[-\mu_0, \mu_0]$.  However, finding an optimal range of $\mu$ seems to be out of reach for our argument, even if we sharpen our estimates, and we therefore prioritize making our argument as simple as possible rather than introducing additional technicalities to squeeze out a marginally more general result.  
\end{REMARK}

\subsubsection{Estimates on $\cF$} 
For distinct $\zeta, \widetilde{\zeta}\in N(\d)$, we have
\begin{align*}
	\frac{|\cF[\zeta](R) - \cF[\widetilde{\zeta}](R)|}{\|\zeta - \widetilde{\zeta}\|_\infty}
	& \le \frac{1}{R^2} \bigg[ \frac{R^2}{3} \|U\|_{L^\infty(I(\d))} + \frac{R^2 \d}{3}\|U'\|_{L^\infty(I(\d))} \bigg] \\
	& \hspace{20 mm} + \bigg[ \frac{4 \|g'''\|_\infty}{g''(1)^2} |V(\l)| + \frac{1}{3g''(1)} \sup_{|\l-1|<\frac{\d}{6}}|V'(\l)|\bigg] \\
	& \le \frac23 \cdot 22(1+M^2) \cdot \frac{10}{g''(1)} + \frac{1}{g''(1)}\cdot (8M + 2) K(\brho, \mu)\\
	& \le \frac{1}{16} + \frac{K(\brho, \mu)}{400}.
\end{align*}

Whenever $|\mu|\le \mu_0$ and $\brho_-(\mu)\le \brho \le \brho_+$, we have $K(\brho, \mu) < 11$, and so the number on the right side in the above inequality less than $9/100$, which in particular is strictly less than $\frac57$.  Since $\|L^{-1}\|=\frac75$, we have in particular that $L^{-1} \cF$ shrinks all distances by a factor that is strictly less than $1$.  

Next, we show that $L^{-1} \cF$ maps $N(\delta)$ to itself.  For $\zeta\in N(\delta)$, we have
\begin{align*}
	\|L^{-1}\cF[\zeta]\|_\infty 
	& \le \|L^{-1} \cF[0]\|_\infty + \|L^{-1}\cF[\zeta] - L^{-1} \cF[0] \|_\infty \\
	& \le \frac35 \frac{K(\brho, \mu)}{g''(1)} + \frac75 \bigg( \frac{1}{16} + \frac{K(\brho, \mu)}{400} \bigg) \|\zeta\|_\infty  \\
	& \le \bigg[ \frac{3}{40} K(\brho, \mu) +  \frac{7}{80} \bigg] \d.
\end{align*}
The number in brackets is $\le 1$ under our assumptions on $\mu$ and $\brho$; we may therefore conclude that $L^{-1} \cF$ maps $N(\d)$ into itself.  We conclude that $L^{-1} \cF$ is a contraction on $N(\d)$, as claimed.  It follows that there exists a fixed point of the map $L^{-1}\cF$, which we denote by $\zeta^{(\brho, \mu)}$.  We let $\f^{(\brho, \mu)}$ denote the corresponding solution of \eqref{e:separated}$_1$, and we also define $\l^{(\brho, \mu)}$ and $y^{(\brho, \mu)}$ from $\f^{(\brho, \mu)}$ in the obvious way.

\subsection{Boundary Conditions}
\label{ss:BC}

Recall that for $|y-1|\le \delta$, we have 
\[
\frac{9}{10} g''(1) \le g''(y) \le \frac{11}{10} g''(1) \qquad \text{ for } |y - 1|<\delta.
\]
Since $g''(1)>0$, The above implies that $g'$ is strictly increasing on $I(\d)$.  Together with the fact that $g'(1) = -\frac13$, it follows that 
\[
\begin{cases}
	g'(y) > 0  & \text{ when } \frac{10}{27g''(1)} = \frac{\d}{27} < y - 1 \le \d, \\
	g'(y) < 0 & \text{ when } -\d\le y-1< \frac{\delta}{33} = \frac{10}{33 g''(1)}.
\end{cases}
\]
Our boundary condition is $g'(y^{(\brho, \mu)}(1)) = 0$.  For fixed $\mu\in [-\mu_0, \mu_0]$, we will show that the numbers $\brho_+$ and $\brho_-(\mu)$ defined in \eqref{e:brho+} and \eqref{e:brho-} satisfy 
\begin{equation} 
	\label{e:ylowerbd}
	y^{(\brho_+, \mu)}(1) - 1 > \frac{\d}{27},
\end{equation} 
\begin{equation} 
	\label{e:yupperbd}
	y^{(\brho_-(\mu), \mu)}(1) - 1 < \frac{\d}{33},
\end{equation} 
and thus that there exists $\brho_0(\mu)\in (\brho_-(\mu), \brho_+)$ such that $g'(y^{(\brho_0(\mu), \mu)}(1)) = 0$.  Our separable solution corresponding to the eigenvalue $\mu$ will thus be $\f^{(\brho_0(\mu), \mu)}$.

\subsubsection{Proof of \eqref{e:ylowerbd}}  
Recall from \eqref{e:epsilonlwr} that $\epsilon(\brho_+, \mu) \ge \frac{19}{2}$. We first show that $\zeta^{(\brho_+, \mu)}$ is bounded below by a positive constant.  We write
\begin{align*}
	\zeta^{(\brho_+, \mu)} & = \frac35 \frac{\epsilon(\brho_+, \mu)}{g''(1)} + L^{-1}\cF[\zeta^{(\brho_+, \mu)}] - L^{-1}\cF[0] \\
	& \ge \bigg( \frac{3}{50} \cdot \frac{19}{2} - \frac75\big( \frac15 + \frac{K(\brho_+, \mu)}{100}\big) \bigg) \delta \\
	& >\frac{\d}{9}. 
\end{align*}
This lower bound on $\zeta^{(\brho_+, \mu)}$ implies a corresponding lower bound on $y^{(\brho_+, \mu)}(1) - 1$.  Indeed, if $\zeta > c$, then 
\begin{align*}
	y(1) - 1 = \int_0^1 \t^2 \zeta(\t)\dtau 
	> \frac{c}{3}.
\end{align*}
In particular, we have the desired inequality on $y^{(\brho_+, \mu)}(1) - 1$, namely
\[
y^{(\brho_+, \mu)}(1) - 1 > \frac{\d}{27}.
\]

\subsubsection{Proof of \eqref{e:yupperbd}}  Recall from \eqref{e:Ksmall} that $K(\brho_-(\mu), \mu)\le \frac{1}{20}$.  Thus

\begin{align*}
	\|\zeta^{(\brho_-(\mu), \mu)}\|_\infty = \|L^{-1}\cF[\zeta^{(\brho_-(\mu), \mu)}]\|_\infty 
	& \le \bigg[ \frac{3}{40} K(\brho_-(\mu), \mu) +  \frac{7}{80} \bigg] \d < \frac{\d}{11}.
\end{align*}

But then 
\begin{align*}
	y^{(\brho_-(\mu), \mu)}(1) - 1 = \int_0^1 \t^2 \zeta^{(\brho_-(\mu), \mu)}(\t)\dtau < \frac{1}{3} \cdot \frac{\d}{11} =  \frac{\d}{33},
\end{align*}
as needed.

\noindent{\bf{Acknowledgments.}}  JJ would like to thank Professor Thomas C. Sideris for bringing her attention to this problem and stimulating discussion. JJ is supported in part by the NSF grant DMS-2306910. TL is supported in part by NSF grant DMS-2408585.  


\end{document}